\documentclass{IEEEtran4PSCC}
\ifCLASSINFOpdf
   \usepackage[pdftex]{graphicx}
\else
   \usepackage[dvips]{graphicx}
\fi
\usepackage[cmex10]{amsmath}
\usepackage{amssymb}
\usepackage{amsfonts}
\usepackage{mathtools}

\usepackage{booktabs}   
\usepackage{multirow}
\usepackage{array}      
\usepackage{arydshln}
\usepackage[numbers]{natbib}
\usepackage{multicol}
\usepackage{algorithm}              
\usepackage[noend]{algpseudocode}   
\usepackage{cite}       
\usepackage{url}        
\usepackage{hyperref}
\usepackage{xcolor}

\makeatletter
\let\old@ps@headings\ps@headings
\let\old@ps@IEEEtitlepagestyle\ps@IEEEtitlepagestyle
\def\psccfooter#1{%
    \def\ps@headings{%
        \old@ps@headings%
        \def\@oddfoot{\strut\hfill#1\hfill\strut}%
        \def\@evenfoot{\strut\hfill#1\hfill\strut}%
    }%
    \def\ps@IEEEtitlepagestyle{%
        \old@ps@IEEEtitlepagestyle%
        \def\@oddfoot{\strut\hfill#1\hfill\strut}%
        \def\@evenfoot{\strut\hfill#1\hfill\strut}%
    }%
    \ps@headings%
}
\makeatother

\psccfooter{%
        \parbox{\textwidth}{\hrulefill \\ \small{24th Power Systems Computation Conference} \hfill \begin{minipage}{0.2\textwidth}\centering \vspace*{4pt} \includegraphics[scale=0.06]{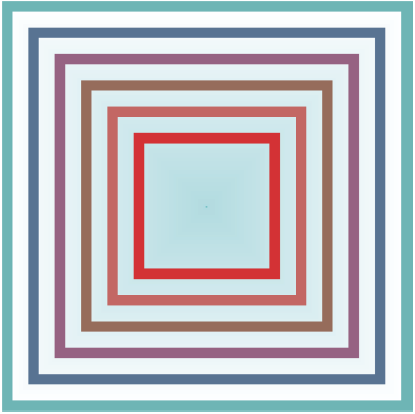}\\\small{PSCC 2026} \end{minipage} \hfill \small{Jingyi Zhao et al. --- June 8-12, 2026}}%
}

\begin{document}
%
\title{Wireless Mobile Charging for Emergency Electric Vehicle
Routing: A Mixed-Integer and Metaheuristic Framework for In-Motion
Energy Transfer}


\author{\IEEEauthorblockN{Jingyi Zhao\IEEEauthorrefmark{1},
Haoxiang Yang\IEEEauthorrefmark{2},
Youxuan Pan\IEEEauthorrefmark{1}, 
 and Yang Liu\IEEEauthorrefmark{3}}
\IEEEauthorblockA{\IEEEauthorrefmark{1}Shenzhen Research Institute of Big Data,  518172, China}
\IEEEauthorblockA{\IEEEauthorrefmark{2}School of Data Science, The Chinese University of Hong Kong, Shenzhen (CUHK-Shenzhen), Shenzhen 518172, China}
\IEEEauthorblockA{\IEEEauthorrefmark{4}Department of Civil and Environmental Engineering,
Department of Industrial Systems Engineering and Management,\\ National University of Singapore, Singapore}
}

\maketitle

\begin{abstract}
As electric vehicles (EVs) become central to decarbonization efforts, the need for uninterrupted power supply in time-critical logistics, particularly in medical transportation, poses unique challenges for power systems integration. Conventional fixed or mobile charging infrastructure requires vehicle downtime, which makes them unsuitable for nonstop operations such as organ delivery. This work introduces the Wireless Mobile Charging Electric Vehicle Routing Problem, a novel framework in which mobile charging trucks wirelessly transfer energy to moving EVs via inductive coupling, eliminating the need for stationary charging stops.
We formulate a mixed-integer programming model that co-optimizes routing and in-motion energy transfer between heterogeneous vehicle fleets under temporal and spatial alignment constraints. To address computational complexity, we develop a hybrid Bitmask Dynamic Programming and Large Neighborhood Search algorithm, tailored to arc-level charging decisions and dual-fleet synchronization. Experiments using real-world hospital logistics data from Singapore demonstrate significant runtime improvements and higher-quality solutions compared with commercial solvers.
This study advances computational power system modeling by incorporating dynamic, motion-based energy delivery that is both time- and space-dependent. The results provide actionable insights into planning mobile energy infrastructure, enhancing demand-side flexibility, strengthening resilience in emergency logistics, and integrating wireless charging into urban power systems to support flexible and zero-carbon electrified mobility.
\end{abstract}



\section{Introduction}
\subsection{Motivation}
The electrification of transport is a cornerstone of decarbonization, but EVs still face limited range and charging delays, restricting their role in both reliable logistics and demand-side flexibility for power systems. Current infrastructure relies on fixed charging stations, which require detours, stationary charging, and impose inflexible load peaks on the grid.  
Mobile charging vehicles (MCVs) offer on-demand solutions: examples include Chargery in Germany, SparkCharge in the U.S. \citep{london_mobile_charger_ev_2024, cnautonews_ev_chargers_2024}, and NIO’s “One-Click Power Up” in China \citep{xinchuxing_nio_charging_service_2024}. While effective in reducing range anxiety, these services still require vehicles to stop and thus remain unsuitable for nonstop operations.  
Recent wireless power transfer prototypes demonstrate in-motion charging through embedded-road systems such as the California pilot \citep{electreon_ers_us_2024}, or vehicle-to-vehicle inductive coupling from Chinese patents and prototypes \citep{CN208164783U, CN111002847A, CN220577090U}. These advances enable continuous charging but have not yet been incorporated into routing and coordination models.  

Existing EVRP studies, e.g., truck–drone coordination \citep{wen2022heterogeneous,mara2023solving} and synchronized mobile charging \citep{Catay2023Matheuristic}, address only static rendezvous. No formulations capture arc-level coordination between heterogeneous fleets in motion. Bridging this gap is essential not only for uninterrupted logistics, but also for quantifying the system-level benefits of wireless mobile charging as a flexible and controllable load resource for future power systems.
\begin{figure}[htbp]
\centering
\begin{minipage}{0.52\textwidth}
    \centering
    \includegraphics[width=\textwidth]{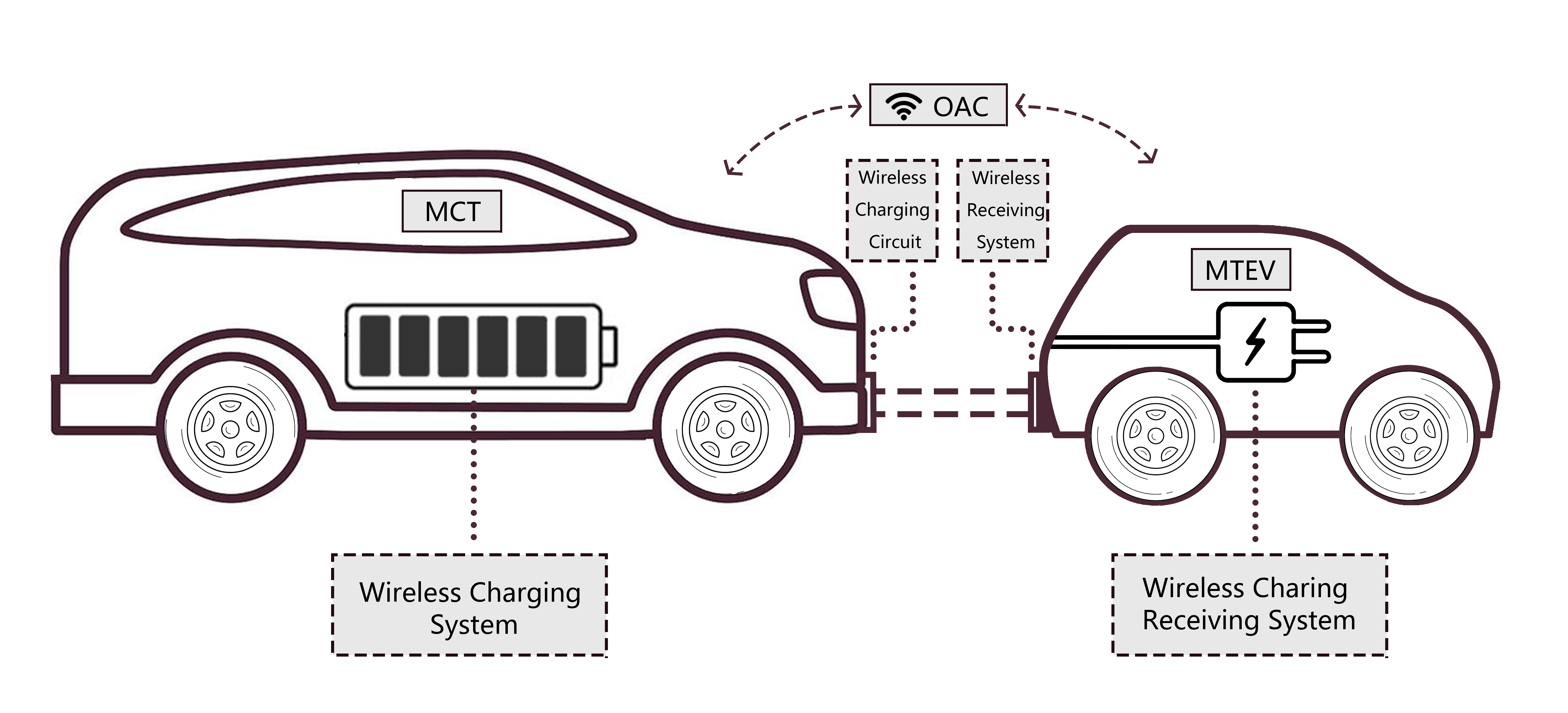}
    \caption{MCT Pattern}
    \label{MCTPattern}
\end{minipage}
\end{figure}

\subsection{Novel Contributions}
This paper makes several contributions at the intersection of electric vehicle routing, wireless energy transfer, and power system flexibility. Compared with existing EV routing formulations that assume charging at static facilities, our work extends the modelling, algorithmic, and system-integration frontier in the following ways:

\begin{enumerate}
    \item \textbf{Problem formulation.} We introduce the \emph{Wireless Mobile Charging Electric Vehicle Routing Problem} (WMC-EVRP), a new class of routing problems where electric delivery vehicles can recharge from mobile charging trucks while both are in motion. To the best of our knowledge, this is the first EVRP variant that models \emph{arc-level} charging decisions coupled with synchronized fleet movement.

    \item \textbf{Mathematical modelling.} We develop a mixed-integer programming (MIP) formulation that jointly optimizes routing and energy transfer between two heterogeneous fleets. The model captures temporal–spatial alignment constraints that are essential for wireless charging during motion, providing a rigorous foundation to evaluate operational feasibility and system benefits.

    \item \textbf{Solution methodology.} To tackle the computational complexity, we design a hybrid algorithm combining Large Neighborhood Search (LNS) with a Bitmask Dynamic Programming (BDP) scheme. While LNS efficiently explores routing neighborhoods, the BDP module resolves binary charging decisions on arcs with high scalability. This integration enables the solution of medium- and large-scale instances that are intractable for commercial solvers.

    \item \textbf{System-level validation.} We perform extensive computational experiments, including realistic urban case studies based on the Singapore road network. Beyond demonstrating computational efficiency, we conduct sensitivity analyses on charging truck operating costs and EV battery capacity. These insights are valuable for both fleet operators and power system planners assessing the integration of wireless mobile charging as a demand-side flexibility option.
\end{enumerate}

Overall, our work establishes the first quantitative framework to assess the operational and system-level impacts of in-motion wireless charging, thereby bridging a critical gap between technological feasibility, optimization modelling, and the needs of future low-carbon electricity systems.

\section{Problem Setting}
We consider the Wireless Mobile Charging Electric Vehicle Routing Problem (WMC-EVRP), defined on a directed graph $\mathcal{G}=(\mathcal{N},\mathcal{A})$ where $\mathcal{N}$ represents the set of depots, hospitals, and potential rendezvous nodes, and $\mathcal{A}$ the travel arcs with associated travel time $t_{ij}$ and \textcolor{black}{distance $c_{ij}$.} Two heterogeneous fleets are modeled:

\begin{itemize}
    \item \textbf{Medical Transport EVs (MTEVs)}: \textcolor{black}{used for urgent deliveries, each with an energy capacity $P$, the maximum load capacity $Q$, and an energy consumption rate $\rho_t$;}
    \item \textbf{Mobile Charging Trucks (MCTs)}: \textcolor{black}{equipped for wireless in-motion charging, each with an energy capacity $B$ and an energy consumption rate $\phi$. When co-traveling with MTEVs, MCTs transfer $\gamma$ as charging efficiency.}
\end{itemize}

The objective is to minimize the total operational cost, including the energy consumption of MTEVs , acquisition cost of MTEVs ($\rho_e$ per deployed vehicle), and acquisition cost of MCTs ($\rho_c$ per deployed vehicle):

\begin{equation}
\min  \sum_{(i,j)\in \mathcal{A}} \sum_{e\in \mathcal{E}}\rho_t c_{ij} x_{ij}^{e}  + \rho_e \sum_{e \in \mathcal{E}} U^e  + \rho_c\sum_{c \in \mathcal{C}} U^c
\end{equation}

The decision variables are defined as follows in Table~\ref{tab:decision_variables}:
\begin{table}[H]
\centering
\caption{Decision Variables}
\begin{tabular}{ll}
\toprule
\textbf{Variable} & \textbf{Description} \\
\midrule
$x_{ij}^{e} \in \{0,1\}$ & 1 if MTEV $e$ travels directly on arc $(i,j)$; 0 otherwise. \\
$z_{ij}^{c} \in \{0,1\}$ & 1 if MCT $c$ travels on arc $(i,j)$; 0 otherwise. \\
$y_i^{e} \in \{0,1\}$ & 1 if MTEV $e$ visits node $i$; 0 otherwise. \\
$w_i^{c} \in \{0,1\}$ & 1 if MCT $c$ visits node $i$; 0 otherwise. \\
$\delta_{ij}^{ec} \in \{0,1\}$ & 1 if MCT $c$ charges MTEV $e$ on arc $(i,j)$; 0 otherwise. \\
$t_i^{e} \in \mathbb{R}_+$ & Arrival time of MTEV $e$ at node $i$. \\
$s_i^{c} \in \mathbb{R}_+$ & Arrival time of MCT $c$ at node $i$. \\
$b_i^{e} \in \mathbb{R}_+$ & Battery level of MTEV $e$ at node $i$. \\
$b_i^{c} \in \mathbb{R}_+$ & Battery level of MCT $c$ at node $i$. \\
$U^e \in \{0,1\}$ & 1 if MTEV $e$ is used (departs from depot); 0 otherwise. \\
$U^c \in \{0,1\}$ & 1 if MCT $c$ is used; 0 otherwise. \\
\bottomrule
\end{tabular}
\label{tab:decision_variables}
\end{table}

The following constraints are imposed, \textcolor{black}{where $M = 10^5$ denotes a sufficiently large constant.}:
\allowdisplaybreaks

\begin{flalign}
& \sum_{j \in \mathcal{N}} x_{ij}^{e} = \sum_{j \in \mathcal{N}} x_{ji}^{e} \qquad \forall i \in \mathcal{N} \setminus \{0, n+1\}, \forall e \in \mathcal{E}. \label{eq1} \\
& \sum_{j \in \mathcal{N}} x_{0j}^{e} = \sum_{j \in \mathcal{N}} x_{j(n+1)}^{e} = U^e \qquad \forall e \in \mathcal{E}. \label{eq2} \\
& \sum_{e \in \mathcal{E}} y_i^e = 1 \qquad \forall i \in \mathcal{N} \setminus \{0, n+1\}. \label{eq3} \\
& \sum_{i\in I \setminus 0}x_{i0}^{e} = \sum_{i\in I \setminus n+1}x_{(n+1)i}^{e} = 0 \qquad \forall e \in \mathcal{E}. \label{eq4} \\
& \sum_{j \in \mathcal{N}} z_{ij}^{c} = \sum_{j \in \mathcal{N}} z_{ji}^{c} \qquad \forall i \in \mathcal{N}, \forall c \in \mathcal{C}. \label{eq5} \\
& \sum_{j \in \mathcal{N}} z_{0j}^{c} = \sum_{j \in \mathcal{N}} z_{j(n+1)}^{c} = U^c \qquad \forall c \in \mathcal{C}. \label{eq6} \\
& t_j^e \geq t_i^e + c_{ij} - M (1 - x_{ij}^e) \qquad \forall (i,j) \in \mathcal{A}, \forall e \in \mathcal{E}. \label{eq7} \\
& s_j^c \geq s_i^c + c_{ij} - M (1 - z_{ij}^c) \qquad \forall (i,j) \in \mathcal{A}, \forall c \in \mathcal{C}. \label{eq8} \\
& s_i^c \leq t_i^e + M (1 - \delta_{ij}^{ec}) \qquad \forall (i,j) \in \mathcal{A}, \forall e \in \mathcal{E}, \forall c \in \mathcal{C}. \label{eq9} \\
& t_0^e = s_0^c = 0 \qquad \forall e \in \mathcal{E}, \forall c \in \mathcal{C}. \label{eq10} \\
& b_j^e = b_i^e - \rho_t c_{ij} x_{ij}^e +\gamma \sum_{c \in \mathcal{C}} \delta_{ij}^{ec}  c_{ij} \quad \forall (i,j) \in \mathcal{A}, \forall e \in \mathcal{E}. \label{eq11} \\  
& b_j^c = b_i^c - \phi c_{ij} z_{ij}^c \qquad \forall (i,j) \in \mathcal{A}, \forall c \in \mathcal{C}. \label{eq12} \\
& 0 \leq b_i^e \leq P \qquad \forall i \in \mathcal{N}, \forall e \in \mathcal{E}. \label{eq13} \\
& \gamma \sum_{j\in I\setminus \{0,i\}}\delta_{ij}^{ec}c_{ij} \leq b_i^c \leq B \quad \forall i \in \mathcal{N}, \forall e \in \mathcal{E}, \forall c \in \mathcal{C}. \label{eq13b} \\
& b_0^e = P,\quad b_0^c = B \qquad \forall e \in \mathcal{E}, \forall c \in \mathcal{C}. \label{eq14} \\
& \sum_{c\in \mathcal{C}}\delta_{ij}^{ec} \leq x_{ij}^{e} \qquad \forall (i,j) \in \mathcal{A}, \forall e \in \mathcal{E}. \label{eq15} \\
& M U^{e} \geq \sum_{i \in I \setminus \{0, n+1\}} y_i^e \qquad \forall e \in \mathcal{E}. \label{eq16} \\
& M U^{c} \geq \sum_{i \in I \setminus \{0, n+1\}} w_i^c \qquad \forall c \in \mathcal{C}. \label{eq16b} \\
& \sum_{i \in \mathcal{N} \setminus \{0, n+1\}} d_i y_i^e \leq Q \qquad \forall e \in \mathcal{E}. \label{eq17}
\end{flalign}

The constraints are interpreted as follows:

Constraint~\eqref{eq1} enforces flow conservation for each MTEV $e$, ensuring that the number of arcs entering and leaving each node is balanced. 
Constraint~\eqref{eq2} guarantees that each MTEV route begins at the departure depot (node $0$) and ends at the return depot (node $n+1$). The binary variable $U^e$ activates the vehicle if it departs the depot.  
Constraint~\eqref{eq3} ensures that each hospital node is visited exactly once by a single MTEV.  
Constraint~\eqref{eq4} prevents MTEVs from returning to node $0$ or starting at node $n+1$, restricting the route to begin and end at the designated depots.  
Constraints~\eqref{eq5} and~\eqref{eq6} impose similar flow and depot usage constraints for MCTs.
Constraints~\eqref{eq7} and~\eqref{eq8} compute the arrival times for MTEVs and MCTs, respectively, based on the travel times and route decisions.  
Constraint~\eqref{eq9} synchronizes charging and thus addresses the coordinate constraints: if MCT $c$ is assigned to charge MTEV $e$ on arc $(i,j)$, then MCT~$c$ must arrive at node $i$ no later than MTEV~$e$.  
Constraint~\eqref{eq10} sets the initial arrival time at the depot to zero.  
Constraint~\eqref{eq11} updates the battery level $b_j^e$ of MTEV $e$ after traversing arc $(i,j)$ by subtracting the energy consumed and adding the energy received from all MCTs charging on that arc.
Similarly, constraint~\eqref{eq12} computes the energy depletion of MCT $c$ after traversing arc $(i,j)$ based on the per-distance consumption rate $\phi$.
Constraints~\eqref{eq13} and ~\eqref{eq13b}  enforces feasible energy states: MTEV and MCT battery levels must stay within their respective capacity at all nodes.
Constraint~\eqref{eq14} initializes all vehicles with fully charged batteries at the depot.  
Constraint~\eqref{eq15} ensures MCT can only provide service on routes traveled by MTEV.
Constrants~\eqref{eq16} and~\eqref{eq16b} activate usage variables when vehicles serve any customer nodes, enabling accurate accounting of fleet costs.  
Finally, constraint~\eqref{eq17} ensures that the total demand served by each MTEV does not exceed its maximum load capacity.

\section{Solution Methodology}

The WMC-EVRP combines two layers of complexity: (i) classical routing constraints of capacitated vehicles under deadlines and (ii) novel arc-level synchronization constraints between MTEVs and MCTs for in-motion wireless charging. Solving the full MIP directly is intractable for realistic city-scale instances, motivating the design of a hybrid algorithm that blends heuristic exploration with exact optimization components. Our approach integrates Large Neighborhood Search (LNS) for routing diversification, Bitmask Dynamic Programming (BDP) for efficient charging-state enumeration, and a reduced Mixed-Integer Program (MIP) for global coordination. This subsection elaborates each component in detail. 

\subsubsection{Large Neighborhood Search for Routing}

LNS forms the backbone of the routing layer. The idea is to iteratively destroy part of the incumbent solution and reconstruct it using specialized heuristics, thereby exploring a large solution space without excessive enumeration.

\paragraph{Removal operators.}
We design six classes of removal  operators tailored for emergency EV routing:
\begin{itemize}
\item \textbf{Random removal (RR):} Randomly selects a group of nodes to remove from the current solution.

\item \textbf{Distance-based removal (DR):}Removes the nodes that add the most distance to the current route.

\item \textbf{String removal (SR):} Removes a segment of nodes from a route. The segment length and starting point are randomly chosen.

\item \textbf{Worst removal (WR):} Removes the nodes that contribute most to the total cost.

\item \textbf{Shaw removal (ShR):} Removes nodes that are similar in terms of location and demand.

\item \textbf{Charge removal (CR):} Removes nodes that cause a vehicle's route to exceed battery capacity by calculating the energy usage and identifying the one with the highest extra consumption.
\end{itemize}

\paragraph{Insertion operators.}
\begin{itemize}
    \item \textbf{Random insertion (RI):} For each removed node, the node is randomly selected and the node is inserted in the position with the lowest cost.

\item \textbf{Greedy insertion (GI):} For each removed node, insert it at the position that causes the smallest increase in cost.

\item \textbf{Sequential insertion (SI):} Inserts each removed node back into the routes one by one.

\item \textbf{Regret-2 insertion (R2I):} For each removed node, compares the best and second-best insertion positions using a regret function.

\item \textbf{Regret-3 insertion (R3I):} Similar to R2I, but also considers the third-best position.

\item \textbf{Charge insertion (CI):} Reinserts removed nodes while ensuring the vehicle’s energy consumption stays within battery limits. \textcolor{black}{This operator can help avoid deploying additional MCTs.}
\end{itemize}

\paragraph{Adaptive mechanism.}
Destroy and repair operators are adaptively weighted based on their historical success rates in improving the incumbent solution, following a roulette-wheel mechanism similar to adaptive LNS in VRP literature.

\subsubsection{Bitmask Dynamic Programming for Arc-level Charging}

Routing decisions alone are insufficient since each feasible MTEV path must be augmented with charging strategies. A naïve enumeration of charging options grows exponentially with the number of arcs. To overcome this, we design a BDP scheme.

\paragraph{State definition.}
We now present the state transition logic used in the BDP algorithm to model battery evolution along a given delivery route. Let $e = 1, \dots, m$ denote the edge index in the current route. At each edge, the algorithm determines whether an MCT charge in transit should occur. The decision for each edge is represented in the $e$-th binary bit of current state $S$, which is set 1 if the MTEV is charge on edge $e$, and 0 otherwise.
Given the battery level $f(e{-}1, S)$ after completing edge $e{-}1$ under state $S$.

\paragraph{Transition.}
The state transitions for edge $e$ follow two possible cases:

\begin{itemize}
    \item \textbf{No Charging on Edge $e$:}  
    The battery level is reduced by the energy consumption $C(e)$:
  $
    f(e, S') = f(e{-}1, S) - C(e),
    $
    where $S' = S$.

    \item \textbf{Charging on Edge $e$:}  
    The battery is increased by $\gamma C(e)$ units (charging gain) while consuming $C(e)$ units of energy to traverse the edge:
  $
    f(e, S'') = \min \bigg( f(e-1, S) + \gamma C(e)-C(e),\ P \bigg),
    $
    where $S''$ is the bitmask obtained by flipping the $e$-th bit of $S$ from $0$ to $1$.
\end{itemize}

\begin{figure}[htbp]
\centering
\includegraphics[width=0.5\textwidth]{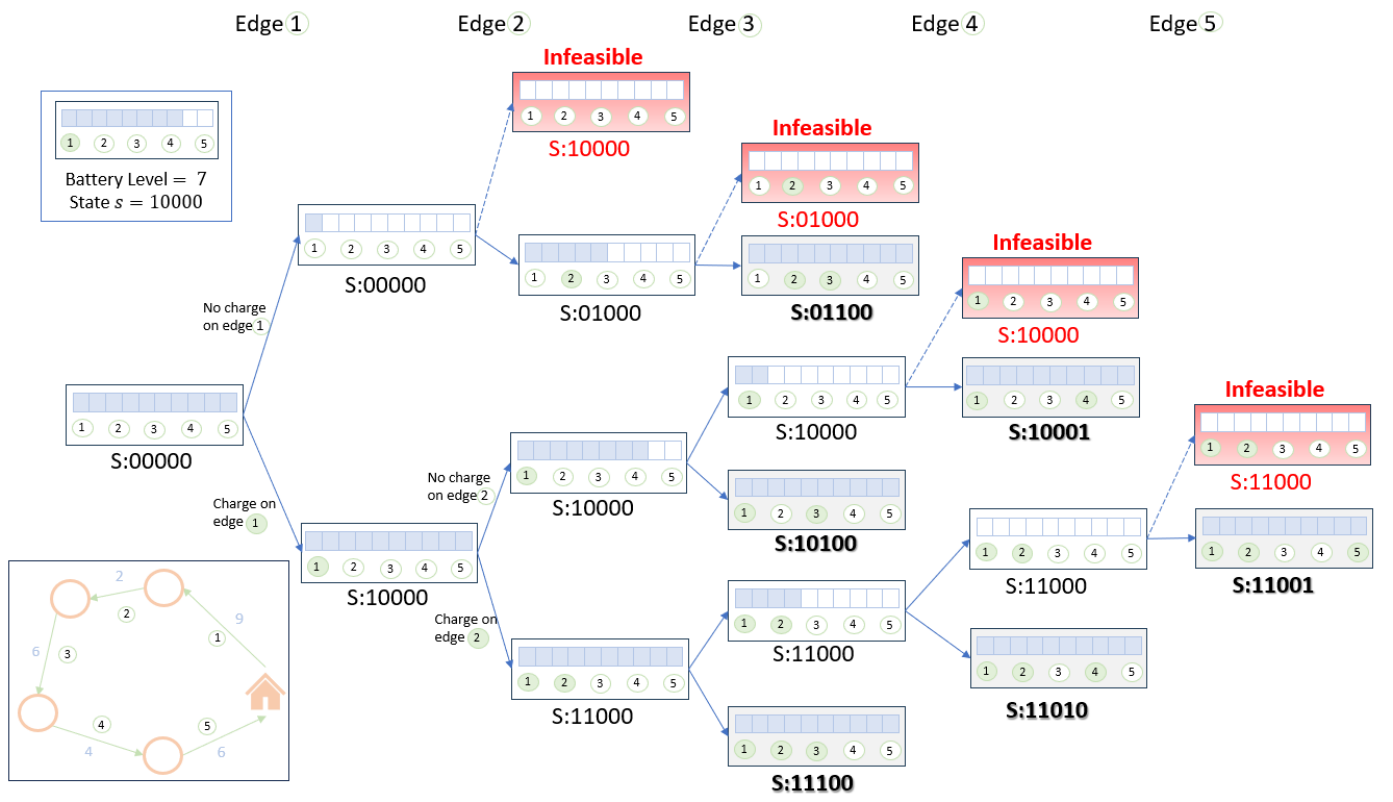}
\caption{The state transition diagram illustrates charging decisions along each edge. Red nodes indicate infeasible states; gray nodes represent terminal feasible states.}
\label{state_trans1}
\end{figure}

Figure~\ref{state_trans1} illustrates the state transitions along a five-edge delivery route under the BDP framework. Each node in the diagram corresponds to a specific charging decision state, encoded as a 5-bit binary string $S$ where each bit indicates whether the MTEV is charged while traversing that edge. The blue horizontal bars next to each state visually represent the remaining battery level of the MTEV after completing the corresponding edge under the current charging state. A solid green border around a bit position indicates that charging occurs on that edge, while an empty border denotes no charging. Nodes highlighted in red represent infeasible states—situations where the MTEV lacks sufficient battery to proceed to the next edge unless charging is performed. These nodes indicate mandatory charging conditions and are pruned from further expansion to prevent energy failure. On the other hand, some states—such as those in the final row of the diagram—do not generate further transitions because the remaining battery at that point is already sufficient to complete the rest of the route. 
These are referred to as \emph{terminal feasible states} and satisfy the condition $f(e, S) \geq r(e)$, \textcolor{black}{where $r(e)$ denotes the minimum remaining battery required to finish the whole route after traversing edge $e$.}
This figure captures the progression of charging decisions, the pruning of infeasible branches, and the early termination of paths when energy sufficiency is reached.

\begin{figure}[htbp]
\centering
\includegraphics[width=0.5\textwidth]{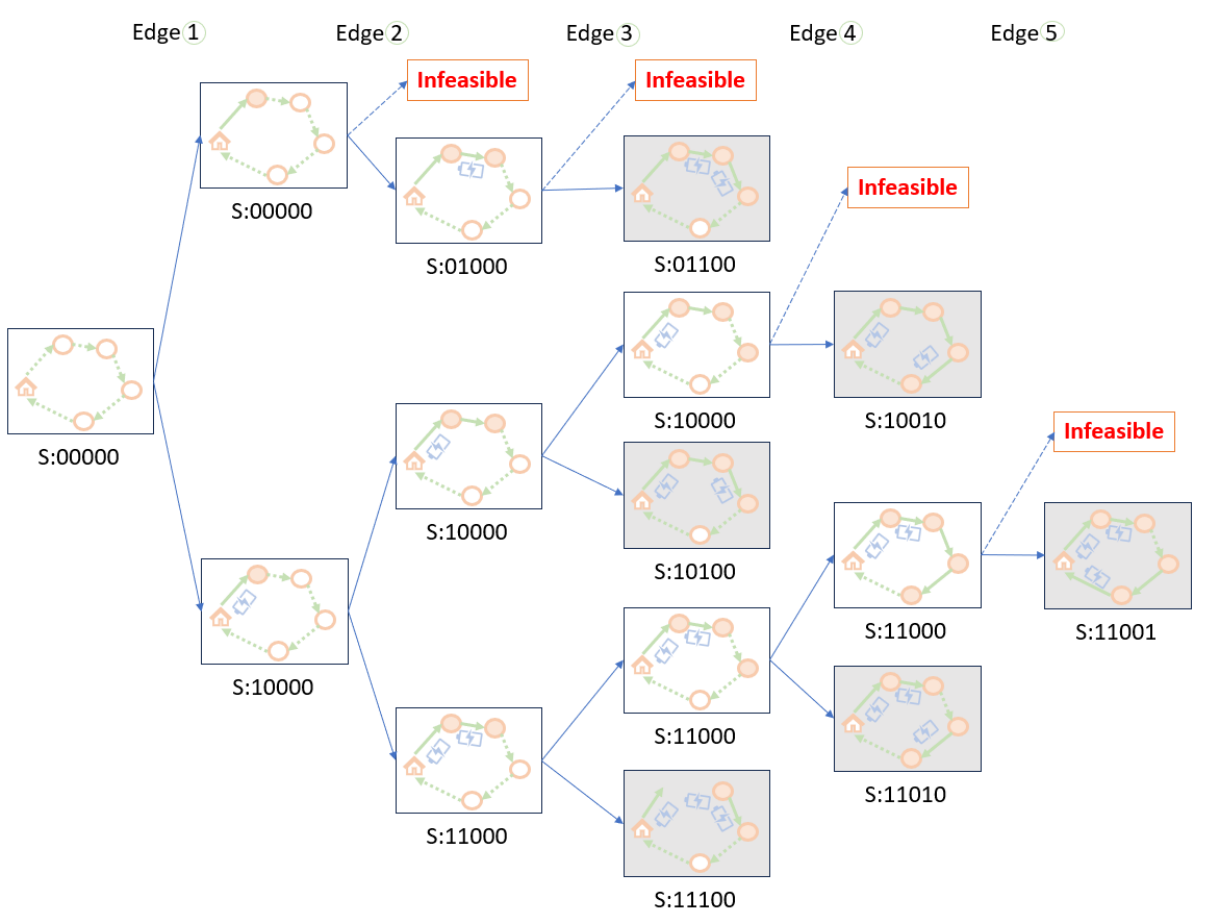}
\caption{Visual example of state transitions for a 5-edge route.}
\label{state_trans2}
\end{figure}

Figure~\ref{state_trans2} illustrates a complete transition tree over a 5-edge route, starting from the initial state $S = \texttt{[00000]}$ with full battery $f(0, S) = P$. Not that for all states, we adopt a left-to-right bit encoding, where the leftmost bit corresponds to edge 1. At each edge, the MTEV can either be charged or not, resulting in two possible branches per state. While the transitions still reflect battery feasibility—as transitions leading to $f(e, S) < 0$ are marked as infeasible and pruned—the primary purpose of this figure is to provide an intuitive visualization of how each bitmask state corresponds to a concrete charging pattern along the route. Each node displays a route diagram with red battery icons on the edges where charging occurs, visually reinforcing the meaning of each binary state. For example, $S = \texttt{[11010]}$  indicates charging on edges 1, 2, and 4. Infeasible states, such as $S = \texttt{[00000]}$ when the battery is insufficient after edge 2, are explicitly labeled and removed. In contrast, the gray-highlighted nodes represent the subset of states that remain feasible and are ultimately retained in the dynamic programming set for further evaluation. These states either reach the end of the route with enough battery or satisfy $f(e, S) \geq r(e)$, thus requiring no additional charging decisions.  

\paragraph{BDP Algorithm Implementation.}
The complete BDP procedure iteratively constructs a battery state table $g(S)$ for each bitmask state $S$, pruning infeasible transitions and retaining only those states that satisfy energy constraints.

\textbf{Space Optimization via Incremental Updates.}
A key observation in the BDP structure is that the battery value of state $S$ at edge $e$ depends only on the value at the previous edge $e{-}1$ under the same bitmask. This reveals a Markovian structure in the transitions, allowing us to reduce the space complexity from $\mathcal{O}(m \cdot 2^m)$ to $\mathcal{O}(2^m)$ , \textcolor{black}{where $m$ is the number of edges in a route. The update proceeds in two steps: first evaluate the battery level assuming charging occurs at edge $e$, followed by the case without charging.}

\textbf{Post-Processing Redundant Solution Removal.}
The in-process marking system described above is effective for pruning redundant solutions with direct parent-child transitions (e.g., $\texttt{[10101]}$ vs. $\texttt{[10100]}$). However, some redundancies exist across independently derived solutions. For example, the two states $\texttt{[10101]}$ and $\texttt{[00101]}$ both represent charging on edges (1, 3, 5) and (3, 5), respectively. Since the latter is a strict subset of the former, the larger state can again be safely pruned.
To remove such cross-solution redundancies, we implement a post-processing step. All feasible states identified during BDP are compared pairwise using a set inclusion test. For any two states $a$ and $b$, if $a \cup b = a$, then $b$ is a subset of $a$, and $a$ is considered redundant. Hence, we can test $  (S_a \,\&\, S_b) = S_a$
to quickly check if $S_a$ is a subset of $S_b$.  In this case, we retain the smaller (minimal) state. 

\textbf{State Marking.}
\textcolor{black}{During the BDP search, each charging configuration is represented by a bitmask state $S$, and its associated battery level is denoted by $g(S)$. To manage feasibility and avoid redundant exploration, we introduce a state marking array $v(S) \in \{-1, 0, 1\}$, which records the status of each configuration based on its battery level:}

\[  
v(S) =  
\begin{cases}  
-1, & \text{if } g(S) < 0 \quad \text{(infeasible; pruned)} \\  
1, & \text{if } g(S) \geq r(e) \quad \text{(terminal feasible)} \\  
0, & \text{otherwise} \quad \text{(active; to be expanded)}  
\end{cases}  
\]

\textcolor{black}{The marking array $v(S)$ allows us to store and manage the classification of each state efficiently during search. Initially, all states are marked as  $v(S) = -1$  (infeasible), except for the root state $S$ where we set $v(S = \texttt{[000...0]})$ to 0. After we explore each new state and compute its battery level $g(S)$, we update $v(S)$ accordingly. This avoids redundant exploration of dominated configurations—for instance, charging on edge 5 when edge 3 already ensures feasibility may lead to an unnecessary state.}

\begin{figure}[htbp]
\centering
\includegraphics[width=0.5\textwidth]{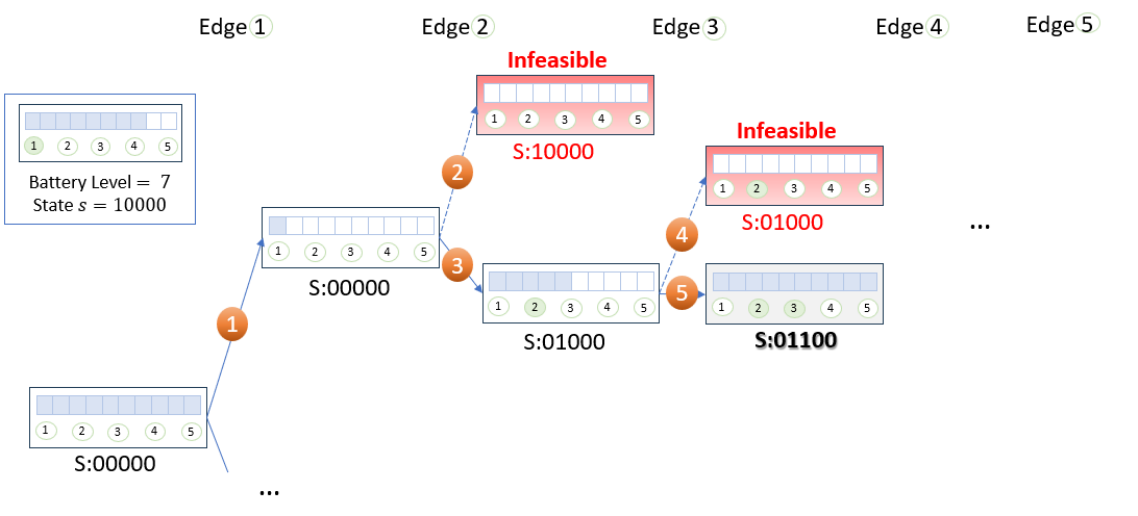}
\caption{
First five BDP transitions. Red denotes infeasible states; numbers indicate update order.
}
\label{state_trans3}
\end{figure}

\textcolor{black}{Figure~\ref{state_trans3} illustrates the update process over a 5-edge route. Below is a step-by-step walkthrough of the first five transitions:}

\begin{enumerate}
    \item \textbf{No-charging on edge 1:}\\
    Compute the battery level after traversing edge 1 without charging:
    \[
        g\bigl(\texttt{[00000]}\bigr) = P - C(1).
    \]
    \item \textbf{No-charging on edge 2:}\\
    Proceed similarly on edge 2:
    \[
        g\bigl(\texttt{[00000]}\bigr) = P - C(1) - C(2).
    \]
    If this value is negative, mark \(v(\texttt{[00000]})\) as infeasible (the vehicle would run out of battery).
    
    \item \textbf{Charging enabled on edge 2:}\\
    Set the bit for edge 2 to 1, indicating in-transit charging:
    \[
        g\bigl(\texttt{[01000]}\bigr) = \min\Bigl(P - C(1) + (\gamma - 1)C(2),\; P\Bigr).
    \]
    This reflects the net battery change after both consuming and gaining energy on edge 2, capped by the maximum battery capacity \(P\).

    \item \textbf{No-charging on edge 3:}\\
    Continue from the state \(\texttt{[01000]}\) without charging on edge 3:
    \[
        g\bigl(\texttt{[01000]}\bigr) = g\bigl(\texttt{[01000]}\bigr) - C(3).
    \]
    If this updated level is negative, mark \(v(\texttt{[01000]})\) as infeasible.

    \item \textbf{Charging on edge 3:}\\
    Flip the bit for edge 3 to indicate charging:
    \[
        g\bigl(\texttt{[01100]}\bigr) = \min\Bigl(g\bigl(\texttt{[01000]}\bigr) + (\gamma - 1)C(3),\; P\Bigr).
    \]
    If \(g\bigl(\texttt{[01100]}\bigr) \ge r(3)\), then \(v(\texttt{[01100]})\) is marked terminal feasible (enough battery remains to finish the route without further charging).
\end{enumerate}

This update process continues until all bitmask states $S$ are classified as either feasible or infeasible. The final output consists of all feasible $S$ satisfying $g(S) \geq 0$, from which optimal charging strategies can be selected. This method exhaustively explores feasible charging states on each route and determines which edges should involve charging, 
ensuring MTEV battery constraints are met.

\textbf{Special Cases in Practice.}
Although the BDP algorithm can theoretically handle routes of any length $m$, 
certain edge cases allow for simple preprocessing checks that reduce unnecessary computation:
\begin{itemize}
    \item If a route has exactly one edge ($m=1$) and its energy consumption $C(0)$ 
    does not exceed the battery capacity $P$, then no in-transit charging is required and the route is trivially feasible. 
    Conversely, if $C(0) > P$, the route is immediately infeasible.
    \item If the total consumption $\sum_{e=0}^{m-1} C(e)$ does not exceed $P$, 
    the vehicle can complete the entire route without any charging. In this scenario, one can bypass the full BDP search 
    and simply record a single no-charging solution.
    \item If some edge $e$ satisfies $C(e) > P$ and even $\gamma \cdot C(e) - C(e)$ is insufficient to keep the battery within its capacity, the route is infeasible from the outset.
\end{itemize}

By performing these checks in constant or linear time before running the BDP, trivial or impossible routes can be pruned, substantially accelerating the overall solution process.

\subsubsection{Reduced MIP for Global Coordination}

After solving the router layer and identifying feasible charging configurations for each MTEV route using BDP, we must select exactly one configuration per route and ensure that these decisions align with MCT movements. This integration step is formulated as a MIP model that builds upon the original vehicle routing constraints ~\eqref{eq1}--\eqref{eq17}.

Let $\mathcal{L}$ be the set of all MTEV routes determined in the routing layer, and $\mathcal{L}_e \in \mathcal{L}$ denote the subset of routes assigned to MTEV $e$. For each route $l\in \mathcal{L}$, the BDP algorithm generates a set $\Omega_l$ of feasible charging bitmask configurations. Each configuration $\omega \in \Omega_l$ represents a specific pattern of charging decisions along route $l$, where each bit indicates whether charging occurs on a particular arc.

We introduce binary variables $Z_{l,\omega} \in \{0,1\}$ to indicate whether configuration $\omega \in \Omega_l$ is selected for route $l$. To establish the connection between the selected charging configurations, MTEV routing, and MCT assignments, we introduce the following constraints:

\begin{align}
    \sum_{\omega \in \Omega_l} Z_{l,\omega} = 1, && \forall l \in \mathcal{L} \label{eq18}
\end{align}
    
\begin{align}
    x_{ij}^e \geq \sum_{l \in \mathcal{L}_e} \sum_{\omega \in \Omega_l : (i,j) \in \omega} Z_{l,\omega} && \forall (i,j) \in \mathcal{A}, \forall e \in \mathcal{E} \label{eq19a}
\end{align}

\begin{align}
    \sum_{c \in \mathcal{C}} \delta_{ij}^{ec} = \sum_{l \in \mathcal{L}_e} \sum_{\omega \in \Omega_l : \omega\text{ charges on }(i,j)} Z_{l,\omega} && \forall (i,j) \in \mathcal{A}, \forall e \in \mathcal{E} \label{eq19b}
\end{align}


Constraint~\eqref{eq18} ensures that exactly one charging configuration is selected for each route. Constraint~\eqref{eq19a} establishes that the routing variables $x_{ij}^{e}$ correctly reflect the arcs used in the selected configurations for MTEV $e$. Constraint~\eqref{eq19b} employs an equality rather than an inequality, ensuring that the sum of charging services provided by all MCTs on arc $(i,j)$ for MTEV $e$ exactly matches the charging requirement specified by the selected configuration. This properly links the BDP decision with the MCT assignments. In addition to these new constraints, the model retains all constraints from the original constraints~\eqref{eq1}--\eqref{eq17}, ensuring feasibility of flow, energy, synchronization, and capacity. This integrated model can be solved using standard MIP solvers such as Gurobi or CPLEX. For large-scale instances where computational complexity becomes prohibitive, we implement a DPS-based heuristic that iteratively selects one configuration per route while systematically pruning infeasible or dominated combinations, significantly reducing the solution space exploration.

\subsubsection{Algorithmic Framework}

The full workflow is summarized in Algorithm~\ref{alg:wmc}, which outlines each step of the LNS-BDP framework.

\begin{algorithm}[htbp]
\caption{Hybrid Algorithm for WMC-EVRP}
\label{alg:wmc}
\begin{algorithmic}[1]
\State Generate initial MTEV routes using nearest-depot heuristic
\Repeat
    \State Select destroy operator and remove subset of customers
    \State Reinsert customers using regret-based repair
    \For{each reconstructed MTEV route}
        \State Run BDP to enumerate feasible charging schedules
    \EndFor
    \State Solve reduced MIP to coordinate MCT deployment
    \State Update incumbent solution if improved
\Until{termination or time limit}
\end{algorithmic}
\end{algorithm}





\section{Numerical Experiments and Validation}

\textcolor{black}{To evaluate the proposed LNS-BDP framework, we conduct a series of numerical experiments. In all settings, the rate $\rho_t$ is set equal to $\phi$, the MCT cost $\rho_c$ is assumed to be twice the MTEV cost $\rho_e$, and the charging efficiency $\gamma$ is set to $2\rho_t$. Each instance is generated by randomly creating a symmetric distance matrix between depot and hospital nodes, and assigning hospital demands randomly between 1 and 3 units. The experiments are structured in three parts: (i) a comparison with Gurobi on small-scale instances, (ii) validation using real-world data from Singapore, and (iii) performance assessment on large-scale instances. All reported results—such as the best ($W_{\mathrm{best}}$) and average ($W_{\mathrm{avg}}$) objective values—are obtained from 10 independent runs of the algorithm. All of the benchmark data, computational results, and detailed route can be downloaded at \url{https://github.com/Jingyi-poly/Coordinated-EV-charging}.}


\subsection{Evaluation of LNS-BDP on Small-Scale Data}

\textcolor{black}{This subsection evaluates algorithm performance on small-scale instances with 5, 9, and 11 customers. Each instance is solved using both our proposed LNS-BDP algorithm and the commercial solver, Gurobi. The results are summarized in Table~\ref{tab:Gurobi}, where Gurobi has a computation time limit of 2~hours per instance.}

\begin{table}[htbp]
\caption{Comparison with Gurobi solver on the small-scaled data}
\label{tab:Gurobi}
\centering
\resizebox{0.5\textwidth}{!}{
\begin{tabular}{ccccccl} 
\toprule
\multirow{2}{*}{\textbf{Instance}} & \multicolumn{2}{c}{\textbf{Gurobi}} & \multicolumn{4}{c}{\textbf{LNS-DP}}                                        \\ 
\cline{2-7}
                                   & O/F & UB                      & $W_\mathrm{best}$ & Gap   & $W_\mathrm{avg}$ &           \\ 
\hline
4A                                 & O   & 6173                    & \textbf{6173}     & 0     & 6173             &                  \\
4B                                 & O   & 4334                    & \textbf{4334}     & 0     & 4334             &                  \\
4C                                 & O   & 6525                    & \textbf{6525}     & 0     & 6525             &                  \\
4D                                 & O   & 5877                    & \textbf{5877}     & 0     & 5883             &                  \\
4E                                 & O   & 5836                    & \textbf{5836}     & 0     & 5878             &                  \\ 
\cdashline{1-6}
8A                                 & O   & 6981                    & \textbf{6981}     & 0     & 6981             &                  \\
8B                                 & O   & 6957                    & \textbf{6957}     & 0     & 6957             &                  \\
8C                                 & O   & 6708                    & \textbf{6708}     & 0     & 6708             &                  \\
8D                                 & O   & 6023                    & \textbf{6023}     & 0     & 6023             &                  \\
8E                                 & O   & 5816                    & \textbf{5816}     & 0     & 5827             &                  \\ 
\cdashline{1-6}
10A                                & O   & 8461                    & \textbf{8461}     & 0     & 8465             &                  \\
10B                                & O   & 7215                    & \textbf{7215}     & 0     & 7218             &                  \\
10C                                & O   & 8633                    & \textbf{8633}     & 0     & 8640             &                  \\
10D                                & F   & 6923                    & \textbf{6923}     & 0     & 6923             &                  \\
10E                                & O   & 7769                    & \textbf{7769}     & 0     & 7802             &                  \\ 
\cdashline{1-6}
12A                                & F   & 7293                    & \textbf{7293}     & 0     & 7300             &                  \\
12B                                & F   & 8651                    & \textbf{8651}     & 0     & 8779             &                  \\
12C                                & F   & 9259                    & \textbf{9224}     & -0.38\% & 9257             &                  \\
12D                                & F   & 8245                    & \textbf{8216}     & -0.35\% & 8218             &                  \\
12E                                & F   & 8373                    & \textbf{8373}     & 0     & 8373             &                  \\
\bottomrule
\end{tabular}
}
\end{table}

\textcolor{black}{Table~\ref{tab:Gurobi} summarizes the computational results. In the table, $O/F$ indicates Gurobi’s solution status—either optimal ($O$) or feasible within the time limit($F$). The column $UB$ denotes the upper bound obtained by Gurobi. The $\mathit{Gap}$ is calculated as $(W_{\mathrm{best}} - \mathrm{UB})/\mathrm{UB} \times 100\%$. As shown in the results, LNS-BDP consistently matches or outperforms the best solutions obtained by Gurobi. Notably, while Gurobi often requires hours to solve an instance, LNS-BDP typically finds comparable solutions within seconds. Moreover, the average objective value $W_{\mathrm{avg}}$ is consistently close to the best value $W_{\mathrm{best}}$, indicating the robustness of the proposed algorithm across multiple runs.}

\subsection{Experiments on Real-World Data}

To test practical relevance, we used hospital location data in Singapore, which features a dense and heterogeneous healthcare network. Singapore’s Human Organ Transplant Act and its advanced medical infrastructure provide an ideal setting for organ transport logistics.  

\begin{table}[htbp]
\centering
\caption{Real-World Hospital Instances in Singapore}
\label{tab:real_world_data}
\resizebox{0.4\textwidth}{!}{
\begin{tabular}{ccccc}
\toprule
\textbf{Instance} & $W_\mathrm{best}$ & $W_\mathrm{avg}$  & $E$ & $C$  \\ 
\midrule
10\_hospital      & 4152              & 4152                      & 2   & 0    \\
11\_hospital      & 3939              & 4050                    & 2   & 0    \\
18\_hospital      & 7215              & 7215                    & 2   & 1    \\
23\_hospital      & 8466              & 8466                    & 3   & 1    \\
26\_hospital      & 8873              & 8880                & 3   & 1    \\
29\_hospital      & 9097              & 9103                  & 3   & 1    \\
\bottomrule
\end{tabular}
}
\end{table}

\textcolor{black}{Table~\ref{tab:real_world_data} presents the results for six real-world instances with 10 to 29 hospitals. For each case, the table reports the total cost, the number of MTEVs ($E$) and MCTs ($C$) deployed. As the number of hospitals increases, both $W_{\mathrm{best}}$ and $W_{\mathrm{avg}}$ grow accordingly, primarily due to longer travel distances. The gap between $W_{\mathrm{avg}}$ and $W_{\mathrm{best}}$ remains consistently small, indicating the robustness of the algorithm across all instances. Meanwhile, we observe that MCTs are gradually introduced as the problem scale increases, reflecting their necessity in maintaining energy feasibility under higher demand.}


\subsection{Experiments on Large-Scale Data}

\textcolor{black}{Finally, we assess scalability and strategic insights on our datasets. Two sets of sensitivity analyses are performed: (i) varying the MTEV battery capacity $P$, and (ii) varying the unit cost of MCTs $\rho_c$. Table~\ref{tab:benchmark} summarizes the average results across three classification types: instance size, $P$, and $\rho_c$. The metrics reported are consistent with those defined in Table~\ref{tab:real_world_data}. Due to space limitations, only aggregated values are shown here; detailed results are available in the supplementary GitHub repository.
} 

\begin{table}[htbp]\caption{Summary Results for LNS-BDP Algorithm Under Variable Parameter Configurations.}
\label{tab:benchmark}
\centering
\resizebox{0.35\textwidth}{!}{
\begin{tabular}{cccccc} 
\toprule
\multicolumn{2}{c}{\textbf{Classification}} & $W_\mathrm{best}$ & $E$ & $C$ \\ 
\cline{1-2}
\textbf{Type} & \textbf{Value} &  &  &  \\ 
\hline
\multirow{6}{*}{Instance} & 60    & 21876  & 6  & 2  \\
                          & 80    & 26895  & 8  & 3  \\
                          & 100   & 32733  & 10 & 3  \\
                          & 120   & 40607  & 13 & 4  \\
                          & 140   & 45441  & 15 & 5  \\
                          & 160   & 54690  & 18 & 7  \\ 
\hline
\multirow{8}{*}{P}        & 400   & 326047 & 69 & 53 \\
                          & 600   & 310579 & 68 & 48 \\
                          & 800   & 279806 & 70 & 37 \\
                          & 1000  & 246719 & 72 & 26 \\
                          & 1200  & 222756 & 70 & 19 \\
                          & 1400  & 190622 & 71 & 8  \\
                          & 1600  & 167461 & 70 & 1  \\ 
\hline
\multirow{7}{*}{$\rho_c$} & 50    & 27472  & 11 & 5  \\
                          & 100   & 27962  & 11 & 4  \\
                          & 500   & 29076  & 12 & 3  \\
                          & 1000  & 30375  & 12 & 3  \\
                          & 1500  & 33065  & 12 & 4  \\
                          & 3000  & 40891  & 12 & 4  \\
                          & 5000  & 43486  & 13 & 3  \\
\bottomrule
\end{tabular}}\end{table}

\begin{figure}[htbp]
    \centering
    \includegraphics[width=0.5\textwidth]{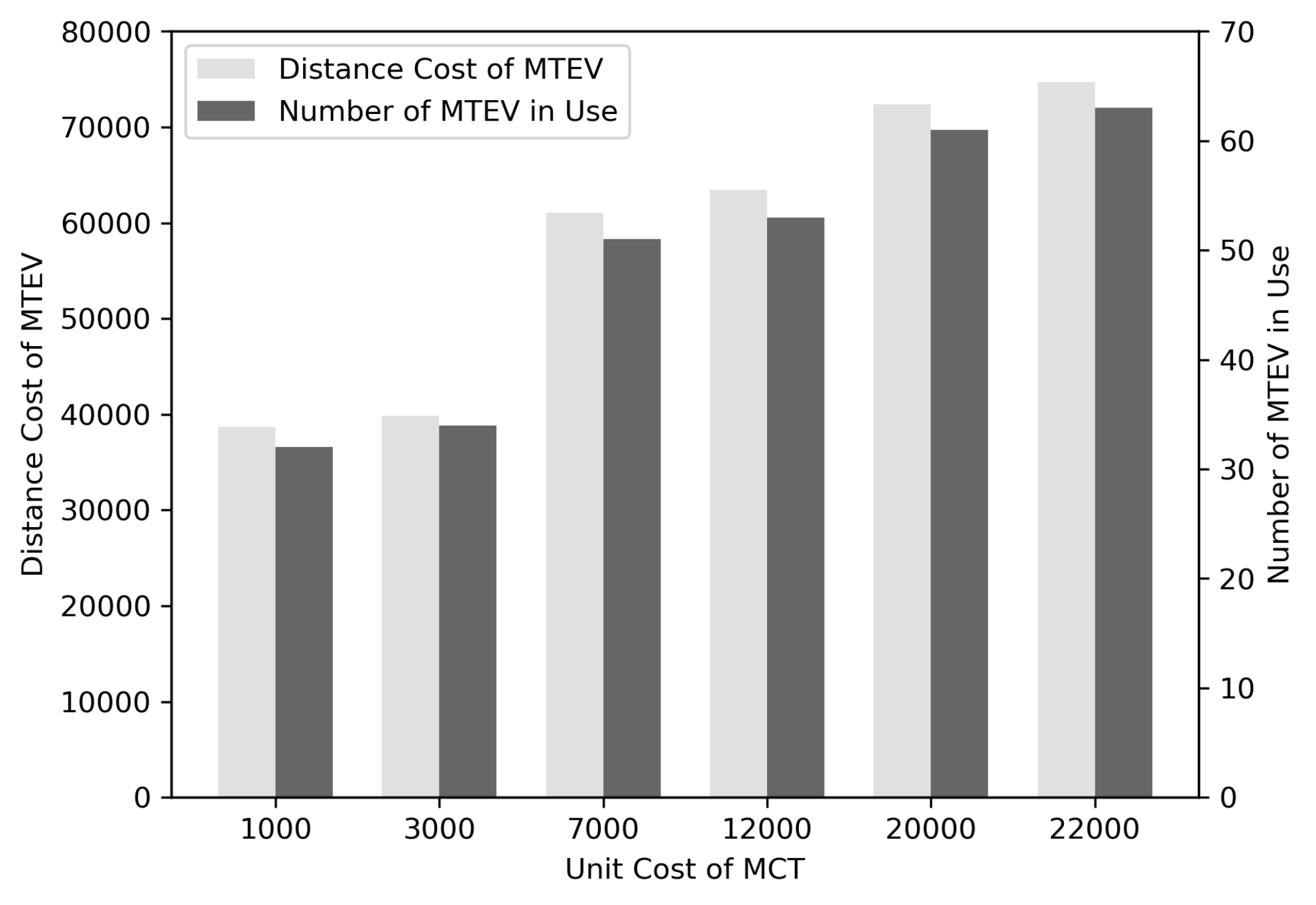}
    \caption{Impact of different unit costs of MCTs ($\rho_c$) on fleet composition and cost.}
    \label{Exp4ChargeCost}
\end{figure}

Figure~\ref{Exp4ChargeCost} illustrates that higher MCT costs push the model toward MTEV-heavy configurations, albeit with higher distance costs.  
Finally, the benchmark analysis in Table~\ref{tab:benchmark} consolidates these findings. \textcolor{black}{As the battery capacity increases from 400 to 1600, the number of required MCTs decreases significantly, while the number of MTEVs remains relatively stable. Specifically, the number of MCTs drops from 53 at ($P = 400$) to just 1 at ($P = 1600$), indicating that larger on-board energy storage substantially reduces the need for en-route charging. In contrast, the number of MTEVs stays within a narrow range of 68–72 across all values of $P$, suggesting that fleet size is primarily driven by demand coverage rather than battery size.}

\section{Conclusion}
This paper proposed the Wireless Mobile Charging Electric Vehicle Routing Problem (WMC-EVRP), a novel model that integrates mobile charging trucks (MCTs) with electric delivery vehicles under in-motion wireless charging. A mixed-integer formulation was developed to capture the temporal–spatial synchronization between heterogeneous fleets, and a hybrid LNS–BDP–MIP algorithm was designed to solve large-scale instances efficiently. Numerical experiments demonstrated substantial improvements in operational efficiency, energy utilization, and cost reduction compared with conventional approaches.  

Beyond logistics applications, the framework highlights the potential of in-motion charging to reshape EV charging demand profiles. By enabling continuous and controllable energy transfer, MCT–MTEV coordination provides new flexibility for power systems, supporting renewable integration, mitigating load peaks, and enhancing system resilience. This study therefore establishes a bridge between emerging wireless charging technologies, advanced optimization models, and the needs of future low-carbon and intelligent power systems.

\textcolor{black}{For future work, incorporating grid-aware and time-dependent elements into the model could enhance its realism and practical value. For example, travel speeds or charging efficiency between nodes may vary with time due to traffic congestion or dynamic grid conditions (e.g., power availability during peak hours). Modeling such spatiotemporal variations would enable more accurate routing and charging strategies that align with both transportation efficiency and power system flexibility.}

\section*{Acknowledgments}
 This work was supported by the Longgang District Special Funds for Science and Technology Innovation under Grant LGKCSDPT2023002.

\end{document}